\documentclass[12pt, reqno]{amsart}
\usepackage[margin=1in]{geometry}

\usepackage{amssymb, mathtools, thmtools, mleftright}
\usepackage[backend=biber, style=numeric, sorting=nyt, giveninits=true, maxbibnames=99]{biblatex}
\usepackage[stretch=10, shrink=10]{microtype}
\usepackage{graphicx}
\usepackage{hyperref}

\DeclareMathOperator{\gl}{GL}

\NewDocumentCommand{\abs}{m}{\mleft\lvert #1 \mright\rvert}
\NewDocumentCommand{\df}{}{\mathop{}\!\mathrm{d}}
\NewDocumentCommand{\dv}{m}{\mathop{}\!\frac{\mathrm{d} #1}{#1}}
\NewDocumentCommand{\numberthis}{}{\addtocounter{equation}{1}\tag{\theequation}}

\RenewDocumentCommand{\Re}{}{\operatorname{Re}}
\RenewDocumentCommand{\Im}{}{\operatorname{Im}}

\numberwithin{equation}{section}
\declaretheorem[numberwithin=section]{theorem}
\declaretheorem[numberlike=theorem]{lemma, corollary, conjecture, proposition}
\theoremstyle{remark}
\newtheorem*{remark}{Remark}

\graphicspath{ {./image} }

\newcommand{\li}{\operatorname{li}}

\begin{document}

\title[On The pair correlation of shifted zeros]
{On The pair correlation of zeros of~\( L \)-functions
  for Non-CM newforms in shifted ranges}

\author{Di Liu}
\address{Di Liu: Department of Mathematics, University of Illinois at Urbana-Champaign,
1409 West Green Street, Urbana, IL 61801, USA}
\email{dil4@illinois.edu}

\author{Clayton Williams}
\address{Clayton Williams: Department of Mathematics, University of Illinois at Urbana-Champaign,
1409 West Green Street, Urbana, IL 61801, USA}
\email{cw78@illinois.edu}

\author{Alexandru Zaharescu}
\address{Alexandru Zaharescu: Department of Mathematics, University of Illinois at Urbana-Champaign,
1409 West Green Street, Urbana, IL 61801, USA \and
Simion Stoilow Institute of Mathematics of the Romanian Academy,
P.O. Box 1-764, RO-014700, Bucharest, Romania}
\email{zaharesc@illinois.edu}

\date{}

\begin{abstract}
  We study the pair correlation between zeros
  of a shifted auxiliary~\( L \)-function attached to a non-CM newform,
  the scale of which is a fixed constant.
  We prove an unconditional asymptotic result for the pair correlation
  and introduce a simplicity hypothesis
  for the zeros of this function, which if true means that multiple zeros of the original
  \(L\)-function cannot be separated by the same fixed distance.
  Our results provide macroscopic information in contrast to
  the pair correlation of the original~\( L \)-function which is of microscopic nature.
\end{abstract}



\maketitle

\section{Introduction}

Let \( \zeta\mleft( s \mright) = \sum_{n = 1}^{\infty} n^{-s} \) for~\( \sigma > 1 \)
denote the Riemann zeta function. The function~\( \zeta\mleft( s \mright) \)
can be analytically continued to the whole complex plane~\( \mathbb{C} \),
except for a simple pole at~\( s = 1 \).
It has trivial zeros at the negative even integers,
and nontrivial zeros inside the critical
strip~\( \mleft\{ s \colon 0 < \Re s < 1 \mright\} \).
The Riemann Hypothesis states that \( \beta = 1/2 \)
for all nontrivial zeros~\( \rho = \beta + i \gamma\).

Assuming this hypothesis,
\citeauthor{montgomeryPairCorrelationZeros1973}~\cite{montgomeryPairCorrelationZeros1973}
studied the pair correlation between
the imaginary parts of the nontrivial zeros. Let \( N\mleft( T \mright) \) be the number of nontrivial
zeros~\( \rho \) with~\( 0 < \gamma < T \); for \( T > 15 \)~and~\( \alpha \in \mathbb{R} \),
define the function
\begin{equation}
  \label{fDef}
  F\mleft( \alpha \mright) \coloneqq \frac{1}{N\mleft( T \mright)}
  \sum_{0 < \gamma, \gamma' < T}
  T^{i \alpha \mleft( \gamma - \gamma' \mright)} w\mleft( \gamma - \gamma' \mright),
\end{equation}
where \( w\mleft( u \mright) \coloneqq 4 / \mleft( 4 + u^{2} \mright) \).
\citeauthor{montgomeryPairCorrelationZeros1973}
showed that when \( 0 \le \alpha < 1 \) we have
\begin{equation}
  \label{eqn:mont-result}
  F\mleft( \alpha \mright)
  = \mleft( 1 + o\mleft( 1 \mright) \mright) T^{- 2 \alpha} \log T
  + \alpha + o\mleft( 1 \mright).
\end{equation}
Later \citeauthor{goldstonLARGEDIFFERENCESCONSECUTIVE1981} extended
the range of~\( \alpha \) to include~\( 1 \)
in~\cite[Lemma~B]{goldstonLARGEDIFFERENCESCONSECUTIVE1981}.
\citeauthor{montgomeryPairCorrelationZeros1973} also conjectured that for \( \alpha \ge 1 \) we have
\begin{equation}
  \label{monConj}
  F\mleft( \alpha \mright) = 1 + o \mleft( 1 \mright).
\end{equation}
Together~\eqref{eqn:mont-result} and the conjectural~\eqref{monConj} completely characterize
the function~\( F\mleft( \alpha \mright) \),
and imply that the normalized imaginary parts
should follow the Gaussian~Unitary~Ensemble (see~\cite[\S 2]{katzZeroesZetaFunctions1999}).
Correlations between tuples beyond pairs, including for~\(L\)-functions beyond the zeta function, have been studied
in~\citeauthor{hejhalTripleCorrelationZeros1994}~\cite{hejhalTripleCorrelationZeros1994},
\citeauthor{rudnickZerosPrincipalLfunctions1996}~\cite{rudnickZerosPrincipalLfunctions1996}, and
\citeauthor{katzZeroesZetaFunctions1999}~\cite{katzZeroesZetaFunctions1999}.

In this paper we consider the pair correlation between zeros of~\( L \)-functions for holomorphic newforms.
Write \( f(z) = \sum_{n > 0} a_{n} n^{\frac{k-1}{2}} q^{n}, q \coloneqq e^{2 \pi i z} \) for \( \Im z > 0 \).
If \( f \) is a cusp form of level~\( N \) with \( a_{1} = 1 \) then \( f \) is a \emph{newform}
if it is a simultaneous eigenform for all Hecke operators, the Fricke involution~\( W_{N} \),
and all Atkin--Lehner involutions~\( W_{Q} \), where \( Q \mid N \)
is a prime~\cite[Definition~2.25]{onoWebModularityArithmetic2004}.
In particular, we are interested in modular forms with trivial nebentypus (see~\autoref{sec: setup}).
The newforms of weight~\( k \) and level~\( N \) form a basis
for the newspace~\( S_{k}^{\text{new}} \mleft( \Gamma_{0}\mleft( N \mright) \mright) \).
Define the \( L \)-function for such a newform, normalized to have symmetry about \( \Re s = 1/2 \),
as the analytic continuation past \( \Re s = 1 \) of the Dirichlet series
\begin{equation*}
  L\mleft( s, f \mright) = \sum_{n \ge 1} \frac{a_n}{n^s}.
\end{equation*}
In particular, \( L\mleft( s, f \mright) \) is a degree~\( 2 \) \( \gl_{2} \)~\( L \)-function,
has an Euler product in the region of absolute convergence,
and corresponds to a completed~\( L \)-function
\begin{equation}
  \label{eqn:completed-L-fcn}
  \Lambda\mleft( s, f \mright) \coloneqq c_k \mleft(2\pi\mright)^{-s}
  \Gamma\mleft( s + \tfrac{k-1}{2} \mright) \mathfrak{q} \mleft(f\mright)^{s/2} L\mleft( s, f \mright),
\end{equation}
which satisfies the functional equation
\begin{equation}
  \label{eqn:fncl-eqn}
  \Lambda\mleft(s, f \mright) = \varepsilon_f \Lambda\mleft( 1-s, f \mright)
\end{equation}
with \( \abs{\varepsilon_f} = 1 \).
Here \( c_k, \mathfrak{q}\mleft(f\mright) \) depend at most on \( N \)~and~\( k \),
and \( \mathfrak{q}\mleft(f\mright) \asymp N k^2 \)~\cite[\S 5.11]{iwaniecAnalyticNumberTheory2004}.

We are interested in correlations between pairs of nontrivial zeros
of~\( L\mleft( s, f \mright) \) whose ordinates are far apart, rather than close.
Similar results for the Riemann zeta function
appear in \citeauthor{chanPairCorrelationZeros2004}~\cite{chanPairCorrelationZeros2004} and
a paper of Ledoan and the third author~\cite{ledoanExplicitFormulasPair2011}.
\citeauthor{chanPairCorrelationZeros2004} uses
analytic tools to redefine \( F\mleft( \alpha \mright) \) on~\cite[182]{chanPairCorrelationZeros2004}
in order to study correlations between distant zeros, while Ledoan and the third author use
number theoretic methods to define an auxiliary~\( L \)-function,
the zeros of which are closely related to those of the Riemann zeta function.
Our work follows the latter paper in its setup.

Throughout this paper, let \( f \in S_k^{\textnormal{new}}\mleft(\Gamma_0\mleft(N\mright)\mright)\)
be a newform without complex multiplication (non-CM, see~\autoref{sec: setup})
and \( L\mleft( s, f \mright) \) be the degree~\( 2 \)~\( L \)-function attached to~\( f \).
For a fixed constant~\( \lambda \ge 0 \), we define a new degree~\( 4 \)~\( L \)-function by
\begin{equation}
  \label{flDef}
  L_{\lambda}\mleft( s, f \mright) \coloneqq
  L\mleft( s + i \tfrac{\lambda}{2}, f \mright)
  L\mleft( s - i \tfrac{\lambda}{2}, f \mright).
\end{equation}
This shifted auxiliary \(L\)-function~\(L_\lambda\mleft(s,f\mright)\)
inherits a functional equation through~\(L\mleft(s,f\mright)\). Write
\begin{equation}
  \Lambda_\lambda\mleft(s,f\mright)
  \coloneqq \Lambda\mleft(s + i \tfrac{\lambda}{2}, f \mright)
  \Lambda\mleft(s - i \tfrac{\lambda}{2}, f\mright);
\end{equation}
then from~\eqref{eqn:fncl-eqn} we have
\begin{equation}
  \label{eqn:shift-fncl-eqn-lambda}
  \Lambda_\lambda\mleft(s,f\mright)
  = \varepsilon_f^2 \Lambda_\lambda\mleft(1-s,f\mright).
\end{equation}

We also need to introduce some auxiliary zero-counting functions.
First define the zero counting functions by
\begin{align*}
  N_{L}\mleft( T \mright)
  & \coloneqq
  \# \mleft\{ \rho \colon L\mleft( \rho, f \mright) = 0,
  0 < \Re \rho < 1, - T < \Im \rho < T \mright\}
  \sim \frac{2}{\pi} T \log T, \\
  N_{L_\lambda}\mleft( T \mright)
  & \coloneqq
  \# \mleft\{ \rho \colon L_\lambda\mleft( \rho, f \mright) = 0,
  0 < \Re \rho < 1, - T < \Im \rho < T \mright\}
  \sim \frac{4}{\pi} T \log T.
\end{align*}
Let \( \rho, \rho' \) run through the nontrivial zeros of~\(L_\lambda\mleft(s,f\mright)\). As in \cite[(2.3)]{murtyExplicitFormulasPair2002} we extend Montgomery's weight function to complex variables by defining\[w(z) := \frac{4}{4-z^2}.\]
For \( \alpha > 0 \) and large~\( T > 0 \) such that \( N_{L_{\lambda}}\mleft( T \mright) > 0 \),
define the shifted pair correlation function
\begin{equation}
  \label{lbdPairFunc}
  F_{\lambda}\mleft( \alpha \mright) = F\mleft( \alpha, L_{\lambda} \mright)
  \coloneqq \frac{1}{N_{L_{\lambda}}\mleft( T \mright)}
  \sum_{-T < \Im \rho, \, \Im \rho' < T}
  T^{4 \alpha\mleft( \rho + \rho' - 1 \mright)} w\mleft( \rho + \rho' - 1 \mright).
\end{equation}

The rationale for introducing this shifted pair correlation function is as follows.
For simplicity we momentarily assume the Riemann hypothesis for~\( L\mleft( s, f \mright) \).
If \( 1 / 2 + i u \)~and~\( 1 / 2 + i u' \) are zeros of~\( \Lambda_\lambda\mleft(s, f\mright) \)
in the critical strip, then \( u = \gamma + \varepsilon  \frac{\lambda}{2} \)~and~\( u' = \gamma' + \varepsilon' \frac{\lambda}{2} \) for some~\( \gamma,\gamma' \)
which are ordinates of zeros of~\( L\mleft( s, f \mright) \)
with~\( \varepsilon, \varepsilon' \in \mleft\{ \pm 1 \mright\} \). By symmetry of the zeros of $L(s,f)$ when $f$ is a newform with trivial nebentypus, we may also write $w(\rho-\rho'-1)$.
Hence the contribution to the weight function
\begin{equation*}
  w\mleft((1/2+iu)-(1/2+iu')-1\mright) = w\mleft(i\mleft(\gamma+\varepsilon\frac{\lambda}{2}\mright)
  -i\mleft(\gamma'+\varepsilon'\frac{\lambda}{2}\mright)\mright)
\end{equation*}
from this pair of zeros is greatest when either \( \varepsilon = \varepsilon' \)
and \( \abs{\gamma-\gamma'} \) is small, or when \( \varepsilon=-\varepsilon' \)
and \( \abs{\gamma-\gamma'} \approx \lambda \).
Hence this shifted pair correlation function measures gaps
on the order of size~\( \lambda \) between zeros of~\( L\mleft(s, f\mright) \).

Our first result is an asymptotic for a sum over the zeros of the shifted auxiliary~\(L\)-function,
where the contribution from those zeros with ordinates separated by~\(\lambda\)
is shown to be nontrivial. Note, in particular,
that our results do not assume the Riemann~Hypothesis
for the function~\( L\mleft( s, f \mright) \).

\begin{remark}
    As a matter of notation, in what follows the constant $\kappa_1$ is $c_{11}>0$, a sufficiently small constant from \cite[Proposition 2.1]{thornerEffectiveFormsSato2021}, when $m=2$ in that proposition. The constant $\kappa_2>0$ depends on $c_{10}$ from \cite[Proposition 2.1]{thornerEffectiveFormsSato2021}, see \eqref{eqn:thetabound}; it also depends on the error term in classic version of the prime number theorem, $\pi(t) = \li(t)+O(te^{-\kappa_2 \sqrt{\log t}})$, \cite[Chapter II.4.1]{tenen-intro}.
\end{remark}

\begin{theorem}
  \label{thm:lbdSum}
  Let \(k, N \in \mathbb{N}\)~with~\(N\) square-free, \(f\) a non-CM holomorphic newform
  with trivial nebentypus in~\(S_k^{\textnormal{new}}\mleft(\Gamma_0\mleft(N\mright)\mright)\),
  \(\lambda>0\), and define \(L_\lambda\mleft(s,f\mright)\) as in~\eqref{flDef}.
  Let \( \rho, \rho' \) run through the nontrivial zeros of~\(L_\lambda\mleft(s,f\mright)\).
  Let \(0 < \delta < 1 / 8\). There exist suitably small constants $\kappa_1,\kappa_2>0$ such that, for \(T \ge 2\), \( T^{\delta} \le x \le T^{1-\delta}\), and \(x\) such that \(2\leq \kappa_1\sqrt{\log x/\log(kN\log x)}\), we have  \begin{multline*}
  \sum_{-T<\operatorname{Im}\rho,\operatorname{Im}\rho'<T}\frac{x^{\rho+\rho'}}{\rho+\rho'} 
= \frac{2}{\pi}\left(1+\frac{\cos(\lambda\log x-\arctan\lambda)}{\sqrt{1+\lambda^2}} \right) Tx\log x \\
\hspace{7em}
-\frac{2}{\pi} \mleft(1-\frac{\left(\lambda ^2-1\right) \cos (\lambda  \log x )-2 \lambda  \sin (\lambda 
   \log x )}{\left(1+\lambda^2\right)^2} \mright)Tx
\\
  + O_{f,\lambda,\delta}\mleft( T x e^{-\kappa_2\sqrt{\log x}}
  \mright).
  \end{multline*}
\end{theorem}

\begin{figure}[htb]
  \centering
  \includegraphics[width=.9\textwidth]{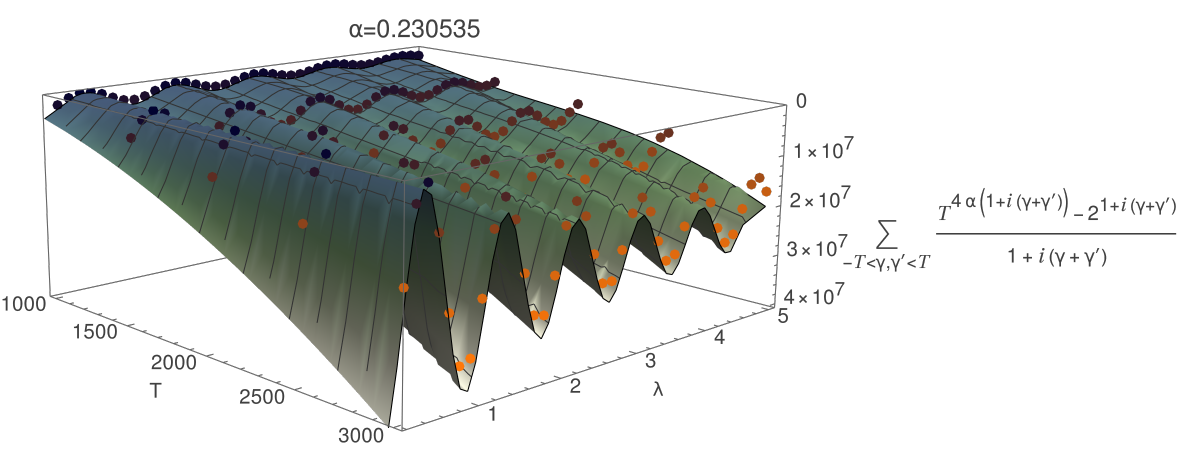}
  \caption{\autoref{thm:lbdSum} in the \( \lambda \)~and~\( T \) aspects for the \( L \)-function
  attached to~\cite[Elliptic~Curve~11.a2]{lmfdb}}
  \label{fig:thm1.1-lambda-aspect}
\end{figure}

In~\autoref{fig:thm1.1-lambda-aspect}, the main terms of~\autoref{thm:lbdSum}
are compared against the sum over zeros for the \( L \)-function attached to
a non-CM elliptic curve~\(E\) for a fixed~\( \alpha \),
where \(E \colon y^2+y=x^3-x^2-10x-20\)~\cite[11.a2]{lmfdb}.
By the modularity theorem such elliptic curves correspond to weight~\( k = 2 \) newforms.

From~\autoref{thm:lbdSum} one can derive an asymptotic for the pair correlation function.

\begin{theorem}
  \label{thm:lbdPair}
  Let \(k, N \in \mathbb{N}\) with \(N\) square-free, \(f\) a non-CM holomorphic newform
  with trivial nebentypus in~\(S_k^{\textnormal{new}}\mleft(\Gamma_0\mleft(N\mright)\mright)\),
  \(\lambda > 0\), \(0 < \delta < 1/8\),
  and define \( F_{\lambda}\mleft( \alpha \mright) \) as in~\eqref{lbdPairFunc}. Let $\kappa_1>0$ be as in \autoref{thm:lbdSum}.  For \( \delta \le \alpha \le 1 / 4 - \delta \)~and~\(T>0\) such that \(N_{L_\lambda}(T)>0\) and \(2\leq \kappa_1\sqrt{4\alpha\log T/\log(4\alpha k N\log T)}\) we have
  \begin{equation*}
    F_{\lambda}\mleft( \alpha \mright)
    = 2\alpha \mleft(1+\frac{4\cos(4\alpha\lambda\log T)}{4+\lambda^2}\mright)+O_{f,\lambda,\delta}\mleft(\frac{1}{\log T}\mright).
  \end{equation*}
\end{theorem}

\begin{figure}[htb]
  \centering
  \includegraphics[width=.7\textwidth]{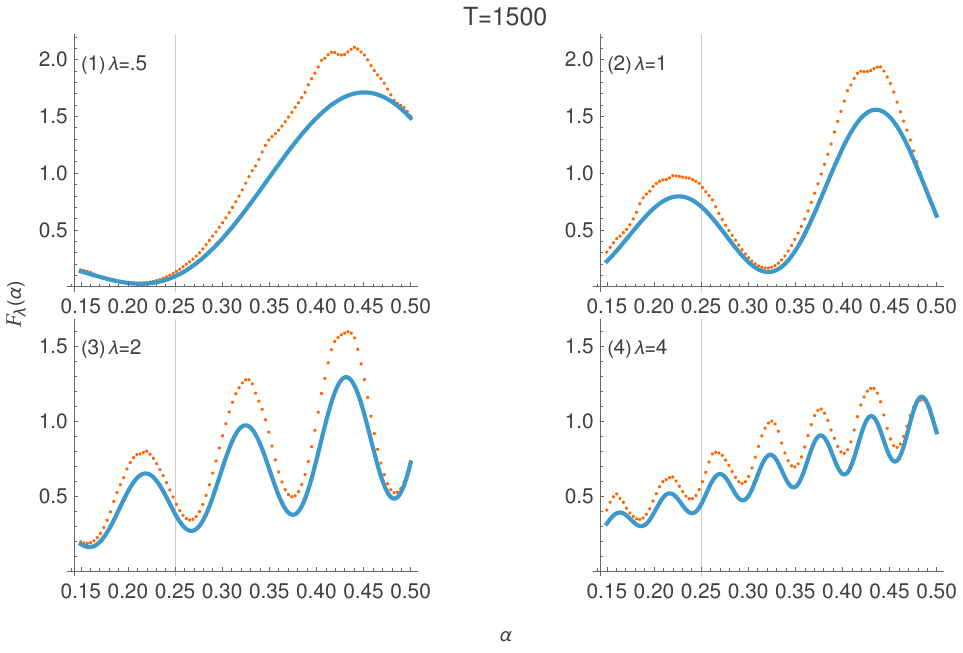}
  \caption{\autoref{thm:lbdPair} in the \( \alpha \)~aspect for small, fixed~\( T \)
  for the \( L \)-function attached to~\cite[Elliptic~Curve~11.a2]{lmfdb}}
  \label{fig:thm1.2-alpha-aspect}
\end{figure}

In~\autoref{fig:thm1.2-alpha-aspect} we plot the asymptotic of~\autoref{thm:lbdPair}
against~\(F_\lambda\mleft(\alpha\mright)\) when~\( \alpha \le 1 / 2 \)
for the same elliptic~curve~\cite[11.a2]{lmfdb}. All figures and data computed with \cite{Mathematica}.

\begin{figure}[htb]
  \centering
  \includegraphics[width=\textwidth]{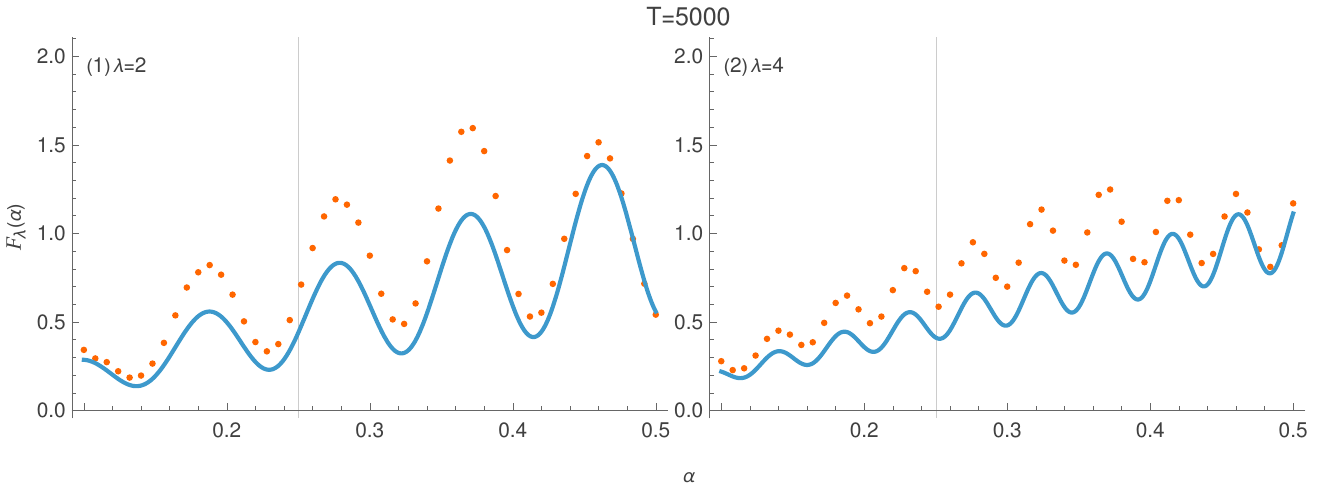}
  \caption{\autoref{thm:lbdPair} for the \(L\)-function attached to~\(\Delta\).}
  \label{fig:Delta-plot}
\end{figure}

An important example of a non-CM newform of level~\( 1 \) and weight~\( k = 12 \)
is the modular discriminant
\begin{equation*}
  \Delta\mleft(z\mright)
  \coloneqq q \prod_{n \ge 1} \mleft(1-q^n\mright)^{24},
\end{equation*}
The \( L \)-function attached to~\(\Delta\)
has the largest ordinate for its lowest zero amongst all primitive algebraic
degree~\( 2 \) \( L \)-functions~\cite[\( L \)-function 2-1-1.1-c11-0-0]{lmfdb}.
We plot~\autoref{thm:lbdPair} for this function in~\autoref{fig:Delta-plot}.

Similar to \citeauthor{rubinsteinChebyshevBias1994}~\cite{rubinsteinChebyshevBias1994},
one can formulate a hypothesis on the simplicity of the zeros of~\(L_\lambda\mleft(s,f\mright)\).
This will in turn be a remark on how many zeros of~\(L\mleft(s,f\mright)\)
can have ordinates separated by the same distance.

\begin{conjecture}
  Assume the Riemann hypothesis for~\(L\mleft(s,f\mright)\) for a given \( f \)
  which is a holomorphic newform. Let \(\lambda>0\).
  Then \(L_\lambda\mleft(s,f\mright)\) has at most one multiple zero with~\(\Im s>\lambda/2\),
  which occurs if and only if there exist two non-trivial zeros~\(\rho, \rho'\)
  of~\(L\mleft(s, f\mright)\) with \(\Im \rho, \Im \rho'>0\) such that
  \begin{equation*}
    \rho-\rho' = i \lambda.
  \end{equation*}
  Moreover, the order of this zero of~\(L_\lambda\mleft(s,f\mright)\),
  if it exists, is \( 2 \).
\end{conjecture}

Because of our normalization the central point is at~\(1/2\) instead of~\(k/2\) for weight~\(k\) newforms.
Note that we avoid the central value~\(s=1/2\) in the conjecture;
this is in contrast to~\cite[176]{rubinsteinChebyshevBias1994} which deals with a different class of~\(L\)-functions.
This is because the central value could be a zero of high multiplicity for~\(L\)-functions of newforms.
For example, in the case of a weight~\( 2 \) newform~\( f \) corresponding to an elliptic curve of large rank,
the Birch~and~Swinnerton-Dyer conjecture~\cite[148]{iwaniecAnalyticNumberTheory2004} would
guarantee \(s=1/2\) is a zero with large order for~\(L\mleft(s, f\mright)\).
It may be the case that elliptic curve~\(L\)-functions
have arbitrarily large orders of zeros at the central point;
\cite{dujellaHistoryEllipticCurves} is a compendium of record-breaking ranks of elliptic curves which, conditional
on the Birch~and~Swinnerton-Dyer conjecture, would give central values with zeros of large orders.
It is known unconditionally that there is an elliptic curve~\(L\)-function
with a zero at~\(s=1/2\) of multiplicity at least~\(3\) in
\citeauthor{grossPointsHeegnerDerivees1983}~\cite{grossPointsHeegnerDerivees1983, grossHeegnerPointsDerivatives1986}
~and~\cite[\S 23.A]{iwaniecAnalyticNumberTheory2004}.

The structure of the paper is as follows. In Section~\ref{sec: setup} we introduce the necessary background,
discussing the pair correlation function and the distribution of coefficients for non-CM newforms.
In Section~\ref{sec:pflg} we introduce and prove a Landau--Gonek formula
for the shifted \( L \)-function. We then use this formula to prove
Theorems~\ref{thm:lbdSum}~and~\ref{thm:lbdPair} in Section ~\ref{sec: main-pfs}.

\section{Setup}
\label{sec: setup}
Recall that a holomorphic modular form~\(f \colon \mathbb{H} \to \mathbb{C}\)
of weight~\(k\) has nebentypus~\(\chi\) a Dirichlet character of modulus~\(N\) if,
for all elements of the group~
\begin{equation*}
  \mleft\{ \begin{pmatrix} a & b \\ c & d \end{pmatrix}
    \colon a, b, c, d \in \mathbb{Z},\,
    N \mid c, \text{ and } ad - bc = 1 \mright\},
\end{equation*}
we have
\begin{equation*}
  f \mleft(\frac{az+b}{cz+d}\mright)
  = \chi\mleft(d\mright) \mleft(cz+d\mright)^k f\mleft(z\mright).
\end{equation*}
If \(\chi\) is the trivial character then one says~\(f\) has trivial nebentypus.

We will need to know the distribution of the coefficients of the newform~\( f \).
We have the following Ramanujan--Petersson type bound
on the Fourier coefficients of such a cusp form from~\citeauthor{deligneConjectureWeil1974}'s
\citeyear{deligneConjectureWeil1974}~proof of the Weil conjectures
~\cite{deligneConjectureWeil1974},

\medskip
\textit{
  Let \( k, N \in \mathbb{N} \) and
  \( \sigma_0\mleft(n\mright) \coloneqq \sum_{d \mid n} 1 \) be the divisor counting function.
  For \( f\mleft(z\mright) = \sum_{n>0} a_n n^{\frac{k-1}{2}} q^n
  \in S_k^{\textnormal{new}}\mleft(\Gamma_0\mleft(N\mright)\mright) \),
  we have~\( \abs{a_n} \le \sigma_0\mleft(n\mright) \).
}
\medskip

When majorizing by primes we are concerned with only the coefficients~\( a_{p} \) for primes~\( p \);
thus we need to know how the~\( a_{p} \)'s distribute in the interval~\( \mleft[ -2, 2 \mright] \).
For each prime~\( p \) define~\( \theta_p \in \mleft[0, \pi\mright] \) by
\begin{equation*}
  a_p = 2 \cos \theta_p.
\end{equation*}
We note that~\(a_n\) is a totally real algebraic integer for all~\(n\),
shown for example in~\cite[Proposition~3.2]{ribetGaloisRepresentationsAttached1977}.
The coefficients with prime index are equidistributed in this interval with respect to one of two measures:
if \( f \) has complex multiplication (CM), then its coefficients are equidistributed with respect to the Hecke distribution,
else they are equidistributed with respect to the Sato--Tate measure.
See also the discussion in~\cite[\S 15.3]{bringmannHarmonicMaassForms2017}.
Recall that a newform~\(f = \sum a_n n^{\frac{k-1}{2}} q^n\) has complex multiplication
by a quadratic character~\(\psi\) of modulus~\(D\)
if \(a_p \psi\mleft(p\mright) = a_p\) for almost all primes~\(p\);
we also say \(f\) has CM by the quadratic field~\( \mathbb{Q}\big( \sqrt{ D } \big) \)
~\cite[Definition~537]{cohenModularForms2017}. In particular,
\( a_p\neq 0\) if and only if~\(p\) totally splits in the ring of integers of the quadratic field.
If no such~\( \psi \) exists then~\( f \) does not have complex multiplication,
and we say~\( f \) is \emph{non-CM}.

For non-CM newforms~\( f \) the angles~\( \theta_{p} \) follow the Sato--Tate distribution,
a well-known problem resolved in~\citeyear{barnet-lambFamilyCalabiYau2011}
by Barnet-Lamb, Geraghty, Harris, and Taylor in~\cite{barnet-lambFamilyCalabiYau2011}.
More recently, in~\citeyear{thornerEffectiveFormsSato2021}
\citeauthor{thornerEffectiveFormsSato2021}~\cite{thornerEffectiveFormsSato2021}
bounded the relative error in the Sato--Tate distribution,
giving the following theorem.

\medskip
\textit{
Let \(N\) be square-free, let \( f \in S_{k}^{\textnormal{new}}\mleft(\Gamma_0\mleft(N\mright)\mright) \)
be a newform with trivial nebentypus and without complex multiplication,
and \(\mleft[\alpha,\beta\mright] \subseteq \mleft[0,\pi\mright]\). Let
\begin{equation*}
  \pi_{f, \mleft[\alpha,\beta\mright]}\mleft(x\mright)
  \coloneqq \#\mleft\{p \le x \colon \theta_p \in \mleft[\alpha,\beta\mright], \, p \nmid N \mright\}
\end{equation*}
and \(\pi\mleft(x\mright)\) be the prime counting function. Then for \(x \ge 3\) we have
\begin{equation}
  \label{sato-tate}
  \frac{\pi_{f,\mleft[\alpha,\beta\mright]}\mleft(x\mright)}{\pi\mleft( x \mright)}
  = \frac{2}{\pi} \int_{\alpha}^{\beta}
  \sin^2\theta \df \theta
  + O\mleft( \frac{\log\mleft(kN \log x\mright)}{\sqrt{\log x}} \mright).
\end{equation}
Here the implied constant is absolute and effectively computable.}
\medskip

We will also need a zero free region for the \( L \)-functions of newforms.
In~\cite[Theorem~5.39]{iwaniecAnalyticNumberTheory2004} we find that
for \(f\) a holomorphic newform of level~\(N\) and weight~\(k \ge 1\),
there exists an absolute constant~\(c>0\)
such that \(L\mleft(s,f\mright)\) has no zero in the region
\begin{equation}
  \label{eqn:zero-free-region}
  \sigma \ge 1-\frac{c}{\log \mleft(N\mleft(\abs{t}+k+3\mright)\mright)},
\end{equation}
except for possibly a single simple real zero~\(\rho = \beta < 1\).

\section{A Landau--Gonek Type Formula}
\label{sec:pflg}

We prove a modified version of
\cite[Lemma~2]{fordDistributionImaginaryParts2009}~and~\cite[Proposition~1]{murtyExplicitFormulasPair2002},
a formula for summing over zeros.
We note that our function~\( \Lambda_\lambda\mleft(s,f\mright) \) satisfies
the hypotheses of~\cite[Lemma~2]{fordDistributionImaginaryParts2009};
namely, that \( L \)-functions of newforms belong to the Selberg class. Denote by \( \Lambda_{L_{\lambda}} \) the Dirichlet coefficients
  of the logarithmic derivative of~\( L_{\lambda}\mleft( s \mright) \) given by
  \begin{equation*}
    \frac{L'_{\lambda}}{L_{\lambda}}\mleft( s \mright)
    = \sum_{n = 1}^{\infty} \frac{\Lambda_{L_{\lambda}}\mleft( n \mright)}{n^{s}},
    \textnormal{ for~} \sigma > 1.
  \end{equation*}
In particular, we have that \( \abs{\Lambda_{L_\lambda}\mleft(n\mright)} \le \Lambda\mleft(n\mright) \ll n^{\varepsilon} \)
for all~\( \varepsilon>0 \), where \( \Lambda\mleft(n\mright) \) is the von-Mangoldt function.
This follows from the proof of the Weil conjectures by
\citeauthor{deligneConjectureWeil1974}~\cite{deligneConjectureWeil1974}.
Taking \( \theta_F \) in~\cite[Lemma~2]{fordDistributionImaginaryParts2009} to be
any~\( 0<\varepsilon<1/2 \) and observing that for
a zero~\( \rho \) of~\( L_\lambda\mleft(s,f\mright)\)
we also have \(L_{\lambda}\mleft( \overline{\rho}, f \mright) = 0 \),
with which one can show the following result, analogous
to~\cite[Proposition~1]{murtyExplicitFormulasPair2002}.

\begin{lemma}
  \label{landaugonek}
  Let \(f\) be a non-CM holomorphic newform
  with trivial nebentypus of integer weight~\(k\),
  \(T\geq 2, x > 1 \) with~\(x \not\in \mathbb{Z}\), and let
  \( L_{\lambda}\mleft( s, f \mright) \) be as in~\eqref{flDef}.
  Let \( \rho = \beta + i \gamma \) run through
  the nontrivial zeros of~\(L_\lambda\mleft(s,f\mright)\).
  Then we have
  \begin{equation*}
    \sum_{-T < \gamma < T} x^{\rho}
    = -\frac{\Lambda_{L_{\lambda}}\mleft( n_{x} \mright)}{\pi}
    \frac{\sin\mleft( T \log \mleft( x / n_{x} \mright) \mright)}
    {\log\mleft( x / n_{x} \mright)}
    + O_{f,\lambda,\varepsilon}\mleft(x^{1+\varepsilon} \log\mleft(Tx\mright)
    +\frac{\log T}{\log x}\mright).
  \end{equation*}
  Here \( n_{x} \) is the closest integer to~\( x \).
\end{lemma}

We briefly sketch how one proves~\autoref{landaugonek}.
We start with the following formula,
obtained by means of the residue theorem:
\begin{equation*}
  \sum_{-T < \gamma < T} x^{\rho}
  = \frac{1}{2 \pi i} \oint_{\mathcal{R}}
  \frac{L_{\lambda}'}{L_{\lambda}}\mleft( s \mright) x^{s} \df s,
\end{equation*}
for some rectangular contour~\( \mathcal{R} \) with
corners \( - \varepsilon / 2 \pm iT \)~and~\( 1 + \varepsilon / 2 \pm iT \),
counterclockwisely labeled \( I_{1}, \ldots, I_{4} \)
starting with \( I_{1} \) as the right vertical contour.

Observe that for large~\( T \) we can deduce from~\cite[(5.27)]{iwaniecAnalyticNumberTheory2004}
that the average spacing of zeros of~\( L_\lambda \mleft( s, f \mright) \)
with~\( \abs{\gamma \pm T} < 1\) is \(\ll 1 / \log \left(\mathfrak{q}\mleft(f\mright)(T+3)^2 \right)\).
Together with~\cite[(5.28)]{iwaniecAnalyticNumberTheory2004},
\begin{equation*}
  \frac{L'_\lambda}{L_\lambda}\mleft(s,f\mright)
  - \sum_{\abs{s - \rho}<1} \frac{1}{s-\rho}
  \ll \log \mleft(\mathfrak{q}\mleft(f\mright) \mleft(T+3\mright)^2\mright),
\end{equation*}
we find that we can always choose a \(T\) such that,
uniformly for all~\( \sigma \in \mleft[-1,2\mright] \), we have
\begin{equation}
  \label{eqn: L'/L bounds horizontal}
  \frac{L_{\lambda}'}{L_{\lambda}}\mleft(\sigma \pm iT\mright) \ll \log^2 T.
\end{equation}
From this we obtain
\begin{equation*}
  I_2, I_4 \ll \mleft(1 + \varepsilon\mright)
  x^{1+\varepsilon/2} \log^2 T.
\end{equation*}

Note that \( \Lambda_{L_{\lambda}}\mleft(n\mright) = 2 \Lambda_L\mleft(n\mright)
\cos\mleft( \tfrac{\lambda}{2}\log n \mright) \).
The main term of~\eqref{landaugonek} comes from the line integral
on~\( \Re s = 1 + \varepsilon / 2 \), whence the logarithmic derivative
can be expanded into its Dirichlet series as
\begin{align*}
  I_{1}
  & = - \frac{1}{2 \pi} \int_{-T}^{T}
  \sum_{n = 1}^{\infty} 2 \Lambda_L\mleft( n \mright)
  \cos\mleft( \tfrac{\lambda}{2} \log n \mright)
  \mleft( \frac{x}{n} \mright)^{1 + \varepsilon / 2 + it} \df t \\
  & = - \frac{1}{\pi} \mleft( \frac{x}{n_{x}} \mright)^{1 + \varepsilon / 2}
  \Lambda_L\mleft( n_x \mright) \cos\mleft( \tfrac{\lambda}{2} \log n_{x} \mright)
  \int_{-T}^{T} \mleft( \frac{x}{n_{x}} \mright)^{it} \df t
  + \mathcal{E},
  \numberthis \label{intOne}
\end{align*}
where the error~\( \mathcal{E} \) is
\begin{equation}
  \label{lgError}
  \mathcal{E} = - \frac{1}{\pi} \sum_{m \ne n_{x}} \int_{-T}^{T}
  \Lambda_L\mleft( m \mright) \cos\mleft( \tfrac{\lambda}{2} \log m \mright)
  \mleft( \frac{x}{m} \mright)^{1 + \varepsilon / 2 + it} \df t.
\end{equation}
The leading term in~\eqref{intOne} gives
\begin{multline*}
  \frac{2}{\pi} \frac{\Lambda_L\mleft( n_{x} \mright) \mleft( x / n_{x} \mright)^{1 + \varepsilon / 2}}
  {\log \mleft( x / n_{x} \mright)} \cos\mleft( \tfrac{\lambda}{2} \log n_{x} \mright)
  \bigg( \frac{\mleft( x / n_{x} \mright)^{iT} - \mleft( x / n_{x} \mright)^{-iT}}{2 i} \bigg) \\
  = \frac{2}{\pi} \frac{\Lambda_L\mleft( n_{x} \mright) \mleft( x / n_{x} \mright)^{1 + \varepsilon / 2}}
  {\log \mleft( x / n_{x} \mright)} \cos\mleft( \tfrac{\lambda}{2} \log n_{x} \mright)
  \sin\mleft( T \log \mleft( x / n_{x} \mright) \mright).
  \numberthis \label{mainOne}
\end{multline*}
Each summand in~\eqref{lgError} can be bounded by
\begin{equation*}
  \abs{ \int_{-T}^{T}
  \Lambda_L\mleft( m \mright) \cos\mleft( \tfrac{\lambda}{2} \log m \mright)
  \mleft( \frac{x}{m} \mright)^{1 + \varepsilon / 2 + it} \df t }
  \le
  2 \abs{ \frac{\Lambda_L\mleft( m \mright) \cos\mleft( \tfrac{\lambda}{2} \log m \mright)}
  {\log \mleft( x / m \mright)} }
  \mleft( \frac{x}{m} \mright)^{1 + \varepsilon / 2}.
\end{equation*}
Proceeding as in \cite[Proposition~1]{murtyExplicitFormulasPair2002}~or~\cite{gonekExplicitFormulaLandau1993}
yields the bound for~\( \mathcal{E} \). To obtain the error corresponding to~\( I_{3} \)
one applies the functional equation~\eqref{eqn:shift-fncl-eqn-lambda} to
the \( L_\lambda' / L_\lambda \)~term in the integrand and proceeds as in~\( I_{1} \),
bounding the contribution from the \( \Gamma \)-functions as in~\cite[Proposition~1]{murtyExplicitFormulasPair2002}.

\section{Proof of~\autoref{thm:lbdSum}}
\label{sec: main-pfs}

We begin with an intermediate proposition.
\begin{proposition}
  Let \( 0 < \varepsilon < 1\)~and~\(2<y<T^{1-\varepsilon}\),
  and assume the hypotheses of~\autoref{landaugonek}. Then
  \begin{multline}
    \label{eqn:prop-eqn-1}
    \abs{ \sum_{-T < \gamma, \gamma' < T}
    \frac{y^{\rho + \rho'} - 2^{\rho + \rho'}}{\rho + \rho'}
    -\int_{2}^{y}
    \frac{\Lambda^{2}_{L_{\lambda}}\mleft( n_{x} \mright)}{\pi^{2}}
    \frac{\sin^{2}\mleft( T \log\mleft( x/n_{x} \mright) \mright)}
    {\log^{2}\mleft( x/n_{x} \mright)} \dv{x}} \\
    \ll_{f, \lambda, \varepsilon} T^{1/2} y^{3/2 + \varepsilon} \log T.
  \end{multline}
\end{proposition}

\begin{proof}
  We begin by writing \( S = \sum_{\rho} x^{\rho},
  M=\frac{\Lambda_{L_\lambda}\mleft(n_x\mright)}{\pi}
  \frac{\sin \mleft(T\log \mleft( x/n_x \mright)\mright)}{\log \mleft( x/n_x \mright)},
  E_1 = x^{1+\varepsilon} \log Tx \)~and~\( E_2=\frac{\log T}{\log x} \).
  Then from~\autoref{landaugonek} we have
  \begin{equation*}
    S^{2} - M^{2} = M\mleft( E_{1} + E_{2} \mright)
    + \mleft( E_{1} + E_{2} \mright)^{2},
  \end{equation*}
  as in~\cite[(4.2)]{murtyExplicitFormulasPair2002}.
  Integrating with respect to~\(x\) from~\(2\) to~\(y\)
  against the kernel function~\(1/x\) yields
  the left hand side of~\eqref{eqn:prop-eqn-1}.
  It remains to bound the error on the right when doing the same operation.
  That is, we aim to bound
  \begin{equation}
    \label{intError}
    \int_2^y \abs{M \mleft( E_1+E_2 \mright)} \dv x
    \text{ and }
    \int_2^y \abs{ E_{1} + E_{2} }^{2} \dv x.
  \end{equation}
  Since \( 2 \abs{ E_{1} E_{2} } \le \abs{ E_{1} }^{2} + \abs{ E_{2} }^{2} \)
  the second integral in~\eqref{intError} is bounded by
  \begin{equation}
    \label{eeBound}
    \int_{2}^{y} x^{1 + 2 \varepsilon} \log^{2}\mleft( T x \mright)
    + \frac{1}{x} \mleft( \frac{\log T}{\log x} \mright)^{2} \df x
    \ll y^{2 + 2 \varepsilon} \log^{2} T.
  \end{equation}

  For the first integral in~\eqref{intError} the Cauchy--Schwarz inequality gives us
  \begin{equation}
    \label{meCS}
    \int_2^y \abs{M \mleft(E_1+E_2\mright)} \dv{x}
    \ll \mleft(\int_2^y \abs{M}^2\dv x\mright)^{1/2}
    \mleft(\int_2^y \abs{ E_1+E_2 }^2 \dv x\mright)^{1/2}.
  \end{equation}
  The integral of~\( \abs{ M }^{2} \) is the same as
  that on the left side of~\eqref{eqn:prop-eqn-1}. 
We have 
\begin{multline*}
  \int_2^y\abs{M}^2 \frac{\df x}{x}= 
\int_2^y \frac{\Lambda_{L_\lambda}^2(n_x)}{\pi^2}\frac{\sin^2 \left(T\log(x/n_x)\right)}{\log^2(x/n_x)}\frac{\df x}{x}
\\
= \sum_{2\leq m\leq y} \frac{\Lambda_\lambda(m)^2}{\pi^2}\int_{m-1/2}^{m+1/2}\frac{\sin^2(T\log(x/m))}{\log^2(x/m)}\frac{\df x}{x}
\end{multline*}
  
  Splitting it into \(O\mleft(\log y\mright)\)
  subintegrals from~\(m-1/2\) to~\(m+1/2, m \in \mathbb{N}\) so that \(n_x=m\)
  on each subinterval yields
  \begin{equation*}
    \frac{\Lambda_{L_\lambda}^2\mleft( m \mright)}{\pi^2}
    \int_{m-1/2}^{m+1/2} \frac{\sin^2\mleft(T\log x/m\mright)}{\log^2\mleft(x/m\mright)} \dv x
    = \frac{T}{\pi} \Lambda_{L_\lambda}^2\mleft(m\mright)
    + O\mleft(m\abs{\Lambda_{L_\lambda}\mleft(m\mright)}^2\mright).
  \end{equation*}

  Summing over~\( m \in \mleft[ 2, y \mright] \) gives
  \begin{equation}
    \label{lambdaSum}
    \int_{2}^{y} \frac{\Lambda^{2}_{L_{\lambda}}\mleft( n_{x} \mright)}{\pi^{2}}
    \frac{\sin^{2}\mleft( T \log\mleft( x/n_{x} \mright) \mright)}
    {\log^{2}\mleft( x/n_{x} \mright)} \dv{x}
    = \frac{T}{\pi} \sum_{2 \le m \le y} \Lambda^{2}_{L_{\lambda}}\mleft( m \mright)
    + O\mleft( y^{2 + \varepsilon} \mright)
    \ll T y^{1 + \varepsilon}.
  \end{equation}
  Now~\eqref{meCS} becomes
  \begin{equation*}
    \int_2^y \abs{M \mleft(E_1+E_2\mright)} \dv{x}
    \ll T^{1/2} y^{1/2 + \varepsilon} \cdot y^{1 + \varepsilon} \log T
    = y^{3/2 + \varepsilon} T^{1/2} \log T,
  \end{equation*}
 which together with~\eqref{eeBound} gives the error in~\eqref{eqn:prop-eqn-1}.
\end{proof}

Next we bound the contribution to~\eqref{eqn:prop-eqn-1}
from \(\sum_{\gamma,\gamma'}\frac{2^{\rho+\rho'}}{\rho+\rho'}\).
Note the summands are large when \( \Im \rho \approx - \Im \rho' \),
so we will split the sum into two cases.
For those with ordinates far apart, we have
\begin{equation}
  \label{twoFarError}
  \abs{ \sum_{\substack{ - T < \Im \rho, \Im \rho' < T \\
    \abs{ \Im \rho - \Im \rho' } > 1 }}
  \frac{2^{\rho + \rho'}}{\rho + \rho'} }
  \le 4 \sum_{- T < \Im \rho < T} \sum_{m = 1}^{2 T} \frac{1}{m}
  \sum_{m < \abs{ \Im \rho' - \Im \rho } \le m + 1} 1
  \ll_{k, N} T \log^{3} T.
\end{equation}
While the pairs with close ordinates are large when \( \Re \rho + \Re \rho' \approx 0 \).
By the functional equation~\eqref{eqn:shift-fncl-eqn-lambda}
and the zero free region~\eqref{eqn:zero-free-region},
\(L_{\lambda}\mleft(s,f\mright)\) has at most two zeros
with real parts close to~\( 0 \), and the width of the zero free region is
\( \min_{\pm} \mleft\{ c / \log \mleft( N \mleft( \abs{ t \pm \lambda } + k + 3 \mright) \mright) \mright\} \),
thus
\begin{equation}
  \label{twoCloseError}
  \abs{ \sum_{\substack{ - T < \Im \rho, \Im \rho' < T \\
    \abs{ \Im \rho - \Im \rho' } \le 1 }}
  \frac{2^{\rho + \rho'}}{\rho + \rho'} }
  \le 4 \sum_{\substack{ - T < \Im \rho < T \\ \abs{ \Im \rho' - \Im \rho } \le 1 }}
  \frac{1}{2 \min \mleft\{ \Re \rho, \Re \rho' \mright\}}
  = O_{k, N}\mleft(T \log^{3} T\mright) + O_f\mleft(1\mright).
\end{equation}

We can evaluate the double sum over the zeros in~\autoref{thm:lbdSum}
now that we have the error from~\eqref{eqn:prop-eqn-1}
as well as the expression in~\eqref{lambdaSum}. Let
\begin{equation*}
  \Psi_{L_{\lambda}}\mleft( x \mright)
  \coloneqq \sum_{n \le x} \Lambda^{2}_{L_{\lambda}}\mleft( n \mright).
\end{equation*}

For our non-CM newform~\( f=\sum a_n n^{\frac{k - 1}{2}} q^n \) write \(a_p = 2\cos\theta_p\).
Using the additive structure of the logarithmic derivatives, one gets
\begin{align*}
  \Psi_{L_{\lambda}}\mleft( x \mright)
  & = \sum_{n \le x} \mleft( \Lambda_{L}\mleft( n \mright)
  \mleft( n^{- i \lambda / 2} + n^{i \lambda / 2} \mright) \mright)^{2} \\
  & = \sum_{p \le x} a_{p}^{2} \cdot \log^{2} p
  \mleft( 2 + 2 \cos\mleft( \lambda \log p \mright) \mright)
  + O\Big( \sum_{\substack{ p^{k} \le x \\ k \ge 2 }} \log^{2} p \Big) \\
  & = \Psi_{1} + \Psi_{2} + O\mleft( \sqrt{ x } \log x \mright),
  \numberthis \label{psiFormula}
\end{align*}
where \( \Psi_{1} = 2 \sum_{p \le x} a_{p}^{2} \log^{2} p \)
~and~\( \Psi_{2} = 2 \sum_{p \le x} a_{p}^{2}
\log^{2} p \cos\mleft( \lambda \log p \mright) \),
which will be evaluated separately using the partial summation formula
\begin{equation}
  \label{partSum}
  \sum_{n \le x} a\mleft( n \mright) f\mleft( n \mright)
  = A\mleft( x \mright) f\mleft( x \mright)
  - \int_{1}^{x} A\mleft( t \mright) f'\mleft( t \mright) \df t,
\end{equation}
with~\( A\mleft( x \mright) \coloneqq \sum_{n \le x} a\mleft( n \mright) \).

For~\( \Psi_{1} \) and \(\Psi_2\), \eqref{partSum} is applied with the summatory function for
\( a\mleft( p \mright) = a_{p}^{2} \).
We rewrite this summatory function~\( A\mleft(x\mright) \) as
\begin{align*}
    A(x) = \pi(x)+\sum_{p\leq x}U_2(\cos\theta_p)+O_f(1),
\end{align*}
where $U_m$ are Chebyshev functions of the second kind, with $U_2(t):=4t^2-1$.  For $f$ a non-CM newform \cite[Proposition 2.1]{thornerEffectiveFormsSato2021} gives an effective bound for the functions
\begin{align*}
    \theta_{f,m}(x):=\sum_{p\leq x}U_m(\cos\theta_p)\log p,
\end{align*}
related to $A(x)$ by partial summation. In our case, Proposition 2.1 (loc. cit.) gives suitably small constants $c_1,c_2,$ and $c_3>0$ such that for all $x$ with $2\leq c_1\sqrt{\log x/\log(k N\log x)}$ we have
\begin{align}\label{eqn:thetabound}
    \abs{\theta_{f,2}(x)}\ll x^{1-c_2} + x \mleft(e^{-c_3\frac{\log x}{\log(2kN)}}+e^{-c_3\sqrt{\log x}}\mright),
\end{align}
from which we derive
\begin{align*}
    \log^2 x\sum_{p\leq x}U_2(\cos\theta_p) = \theta_{f,2}(x)\log x-\log^2 x\int_2^x\frac{\theta_{f,2}(t)}{t\log^2 t}\df t\ll xe^{-c_3\sqrt{\log x}}
\end{align*}
for suitably large $x>0$.

Let $c_4>0$ come from the error term in the classic form of the prime number theorem, so $\pi(t) = \li(t)+O(te^{-c_4\sqrt{\log t}}) $.  Hence for $\Psi_1$ with $A(x)$ as above and $f(t) = 2\log^2 t$, we get
\begin{align*}
    \Psi_1(x) &= 2\mleft(\pi(x)+\sum_{p\leq x}U_2(\cos\theta_p)\mright)\log^2 x-4\int_1^x\frac{\mleft(\pi(t)+\sum_{p\leq x}U_2(\cos\theta_p)\mright)\log t}{t}\df t.
\end{align*}
Because $\int_1^x \frac{\pi(t)\log t}{t}\df t = \frac{1}{2}(x-x\log x+\li(x)\log^2 x-1)+O(xe^{-c_4\sqrt{\log x}})$, there is a $c_5>0$ such that
\begin{align}\label{eqn:Psi1-error}
   \Psi_1(x) = 2x\log x-2x+O_f(xe^{-c_5\sqrt{\log x}})
\end{align}
for sufficiently large $x>0$.\par

Our analysis for $\psi_2$ is similar. Substituting our expression for $A$ into \eqref{partSum} with $f(t) = 2\log^2 t\cos t(\lambda\log t)$ we obtain
\begin{multline*}
\Psi_2(x)=2 \li x \log^2 x \,\cos(\lambda\log x)
-2\int_{1}^x \frac{\li(t)\log t}{t}\left(2\cos(\lambda\log t)-\lambda \log t\, \sin(\lambda\log t)\right) \df t\\
+O_{f,\lambda}\left(xe^{-c_5\sqrt{\log x}}\right).
\end{multline*}
It remains to evaluate the integral. Substitution yields
\begin{align*}
\int_{\log 2}^{\log x}2\li \left(e^u\right) u \cos(\lambda u)-\lambda \li \left(e^u\right)u^2 \sin(\lambda u)\df u.
\end{align*}
Now $\operatorname{li}(e^u)$ is related to the exponential integral by \cite[p. xxxiii]{gradshteinTableIntegralsSeries2015}.
\begin{align*}
\li \left(e^u\right) = \operatorname{Ei}(u):=-\int_{-u}^\infty e^{-t}t^{-1}\df t.
\end{align*}

Evaluating the integrals symbolically with Mathematica \cite{Mathematica} gives the following. Intermediate terms containing $\operatorname{Ei}$ cancel in the final evaluation.

\begin{multline*}
\int_{0}^{\log x}2\operatorname{Ei} \left(u\right) u \cos(\lambda u)-\lambda \operatorname{Ei} \left(u\right)u^2 \sin(\lambda u)\df u\\
=\li(x) \log ^2(x) \cos (\lambda  \log x)
-
\frac{x \log x}{1+\lambda^2} \left(\lambda  \sin (\lambda  \log x)+\cos (\lambda 
   \log x)\right)\\
   +x\Bigg(\frac{2 \lambda  \sin (\lambda  \log
   (x))}{\left(1+\lambda^2\right)^2}+\frac{\left(1-\lambda ^2\right) \cos (\lambda  \log x)}{\left(1+\lambda
   ^2\right)^2}\Bigg)+O_\lambda(1).
\end{multline*}
Hence for sufficiently large $x$ we have
\begin{multline}
\Psi_2(x) = \frac{2 x\log x}{1+\lambda^2} \Big(\lambda  \sin (\lambda  \log
   (x))+\cos (\lambda  \log x)\Big)\\
   +\frac{2x}{\left(1+\lambda^2\right)^2} \Big( \left(\lambda ^2-1\right) \cos (\lambda  \log x)-2 \lambda  \sin (\lambda 
   \log x)\Big)\\\
   +O_\lambda(xe^{-c_6\sqrt{\log x}}).\label{eqn:Psi2-error}
\end{multline} 

Equations \eqref{eqn:Psi1-error} and \eqref{eqn:Psi2-error}, together with the errors~\eqref{twoFarError},
\eqref{twoCloseError}~and~\eqref{eqn:prop-eqn-1}, yields for sufficiently large $T$ and $x$ the asymptotic 
\begin{multline}
  \label{sumAsymp}
\sum_{-T<\operatorname{Im}\rho,\operatorname{Im}\rho'<T}\frac{x^{\rho+\rho'}}{\rho+\rho'} 
= \frac{2}{\pi}\left(1+\frac{\cos(\lambda\log x-\arctan\lambda)}{\sqrt{1+\lambda^2}} \right) Tx\log x \\
-\frac{2}{\pi} \mleft(1-\frac{1}{\left(1+\lambda^2\right)^2} \Big( \left(\lambda ^2-1\right) \cos (\lambda  \log x )-2 \lambda  \sin (\lambda 
   \log x )\Big)\mright)Tx
\\
  + O_{f,\varepsilon}\mleft( T x e^{-c_7\sqrt{\log x}}
  + T^{1/2} x^{3/2+\varepsilon} \log T +T \log^3 T\mright).
\end{multline}
Note the first term inside~\(O_{f,\varepsilon}\) dominates the other two
when~\( T^{\delta} \le x \le T^{1 - \delta} \),
and this concludes the proof of~\autoref{thm:lbdSum}.

\section{Proof of~\autoref{thm:lbdPair}}
\label{sec:main-pfs-1}

Now we proceed to include the weight function~\( w\mleft( iu \mright) \),
following the method in~\cite[Section~6]{murtyExplicitFormulasPair2002}. 
Let
\begin{equation*}
  h_{x}\mleft( y \mright) \coloneqq
  \begin{cases}
    - 1 / x, & \text{ if } y \le x, \\
    3 x^{3} / y^{4}, & \text{ if } y > x.
  \end{cases}
\end{equation*}
After changing the variable in~\eqref{sumAsymp} from \( x \) to~\( y \),
multiplying by~\( h_{x}\mleft( y \mright) \) and integrating,
each summand on the left hand side becomes
\begin{align}
  \frac{1}{\rho + \rho'} \mleft( - \frac{1}{x} \int_{1}^{x} y^{\rho + \rho'} \df y
  + 3 x^{3} \int_{x}^{\infty} y^{\rho + \rho' - 4} \df y \mright)
  & = - \frac{x^{\rho + \rho'}}{\rho + \rho'}
  \mleft( \frac{1}{\rho + \rho' + 1} + \frac{3}{\rho + \rho' - 3} \mright)
  + \mathcal{E} \\
  & = - \frac{4 x^{\rho + \rho'}}{\mleft( \rho + \rho' + 1 \mright)
  \mleft( \rho + \rho' - 3 \mright)} + \mathcal{E} \\
  & = x^{\rho + \rho'} w\mleft( \rho + \rho' - 1 \mright) + \mathcal{E},
  \label{eqn:lhshxy}
\end{align}
where \( \mathcal{E} = \mleft( x \mleft( \rho + \rho' \mright) \mleft( \rho + \rho' + 1 \mright) \mright)^{-1} \)
and satisfies
\begin{equation*}
  \sum_{- T < \Im \rho, \, \Im \rho' < T} \mathcal{E}
  \ll_{f} \frac{T}{x} \log^{2} T.
\end{equation*}

Next, we rewrite the main term on the right hand side of~\eqref{sumAsymp} as
\begin{multline}
  \label{phiPart}
 \frac{2}{\pi}T\Bigg[y\log y+\frac{y\log y\cos(\lambda\log y-\arctan\lambda)}{\sqrt{1+\lambda^2}} -y
 \\+\frac{y}{\left(1+\lambda^2\right)^2} \Big( \left(\lambda ^2-1\right) \cos (\lambda  \log y )-2 \lambda  \sin (\lambda 
   \log y )\Big)\Bigg]\\
  = \frac{2}{\pi} T \mleft( \Phi_{1}\mleft( y \mright)
  + \frac{1}{\sqrt{ 1 + \lambda^{2} }} \Phi_{2}\mleft( y \mright)+\Phi_3(y)+\frac{1}{(1+\lambda^2)^2}\Phi_4(y) \mright)\\
  =\frac{2}{\pi} T \mleft( \Phi_{1}\mleft( y \mright)
  + \frac{1}{\sqrt{ 1 + \lambda^{2} }} \Phi_{2}\mleft( y \mright) \mright)+O_\lambda(y).
\end{multline}
We first remark that \( \int \Phi_{j}(y) h_{x}(y)\df y \ll_\lambda x\) for \(j=3,4\), and we discard these terms. 
The integral~\( \int \Phi_{1}(y) h_{x}(y)\df y \) is
\begin{equation}
  \label{leadAsymp}
  \begin{split}
   - \frac{1}{x}\int_{1}^{x}y \log y \df y + 3 x^{3} \int_{x}^{\infty} \frac{\log y}{y^{3}} \df y
    & =-\mleft( \frac{1}{2}x\log x-\frac{1}{4}x+\frac{1}{4x}\mright)+\mleft(\frac{3}{2}x\log x+\frac{3}{4}x\mright)
    \\
    & = x\log x+x-\frac{1}{4x}.
  \end{split}
\end{equation}
Write $\theta=\arctan\lambda$.
On the other hand, the integral~\( \int \Phi_{2}(y) h_{x}(y)\df y \) is
\begin{equation}
  \label{cosInt}
  -\frac{1}{x} \int_{1}^{x} y\log y \cos\mleft( \lambda \log y - \theta \mright) \df y
  + 3 x^{3} \int_{x}^{\infty} \frac{\log y}{y^{3}}
  \cos\mleft( \lambda \log y - \theta \mright) \df y
  = \Phi_{2, 1} +  \Phi_{2, 2}.
\end{equation}
The first integral in~\eqref{cosInt} is
\begin{multline*}
  \Phi_{2, 1}
  = -\frac{1}{x} \cos \theta \int_{1}^{x} y\log y \cos\mleft( \lambda \log y \mright) \df y
  -\frac{1}{x} \sin \theta \int_{1}^{x} y \log y \sin\mleft( \lambda \log y \mright) \df y \\
  = -\frac{1}{x} \cos \theta \int_{0}^{\log x} u e^{2u} \cos\mleft( \lambda u \mright) \df u-
  \frac{1}{x} \sin \theta \int_{0}^{\log x} u e^{2u} \sin\mleft( \lambda u \mright) \df u.
\end{multline*}
By~\cite[(2.667.5) \& (2.667.6)]{gradshteinTableIntegralsSeries2015} this is
\begin{align}
  \Phi_{2, 1} &= \mleft. - \frac{e^{2u}}{4 + \lambda^{2}}
  \mleft( \mleft( 2u - \frac{4 - \lambda^{2}}{4 + \lambda^{2}} \mright)
  \cos\mleft( \lambda u - \theta \mright)
  + \mleft( \lambda u - \frac{4 \lambda}{4 + \lambda^{2}} \mright)
  \sin \mleft( \lambda u - \theta \mright) \mright) \mright|_{0}^{\log x}\nonumber\\
  &=
-x \log x\frac{ (2 \cos (\lambda\log x-\theta)+\lambda  \sin (\lambda\log x-\theta))}{4+\lambda^2}+O_\lambda(x).
  \numberthis \label{cos1st}
\end{align}
The second integral in~\eqref{cosInt} is
\begin{multline*}
  \Phi_{2, 2}
  = 3 x^{3}\cos \theta \int_{x}^{\infty} \frac{\log y}{y^{3}}
  \cos\mleft( \lambda \log y \mright) \df y
  + 3 x^{3}\sin \theta \int_{x}^{\infty} \frac{\log y}{y^{3}}
  \sin\mleft( \lambda \log y \mright) \df y \\
  =3 x^{3} \cos \theta \int_{\log x}^{\infty} u e^{-2 u}
  \cos\mleft( \lambda u \mright) \df u
  + 3 x^{3}\sin \theta \int_{\log x}^{\infty} u e^{-2 u}
  \sin\mleft( \lambda u \mright) \df u.
\end{multline*}
Again, by~\cite[(2.667.5) \& (2.667.6)]{gradshteinTableIntegralsSeries2015}, this equals
\begin{multline*}
  \mleft. 3 x^{3} \cdot \frac{e^{-2 u}}{4 + \lambda^{2}}
  \mleft( \mleft( -2 u - \frac{4 - \lambda^{2}}{4 + \lambda^{2}} \mright)
  \cos\mleft( \lambda u - \theta \mright)
  + \mleft( \lambda u + \frac{4 \lambda}{4 + \lambda^{2}} \mright)
  \sin\mleft( \lambda u - \theta \mright) \mright) \mright|_{\log x}^{\infty} \\
  = \frac{3 x}{4 + \lambda^{2}}
  \mleft( \mleft( 2 \log x + \frac{4 - \lambda^{2}}{4 + \lambda^{2}} \mright)
  \cos\mleft( \lambda \log x - \theta \mright)
  - \mleft( \lambda \log x + \frac{4 \lambda}{4 + \lambda^{2}} \mright)
  \sin\mleft( \lambda \log x - \theta \mright) \mright).
  \numberthis \label{cos2nd}
\end{multline*}

The largest terms in \eqref{cos1st}~and~\eqref{cos2nd}
are of size~\( x \log x \), which give
\begin{align*}
  %\label{cosAsymp}
  \Phi_{2, 1} +  \Phi_{2, 2}
  &= \frac{4x \log x}{4+\lambda^2} \mleft(
  \cos(\lambda\log x-\theta)-\lambda\sin(\lambda\log x-\theta)
  \mright)+O_\lambda(x).
  \end{align*}
We get
\begin{align}\label{eqn:intrhshxy}
  \frac{2T}{\pi}\int_1^\infty\Bigg(\Phi_1(y)&+\frac{1}{\sqrt{1+\lambda^2}}\Phi_2(y)\Bigg)h_x(y)\df y
  \nonumber\\
  &=\frac{2Tx\log x}{\pi} \mleft(
  1+\frac{4\cos(\lambda\log x)}{(4+\lambda^2)}
  \mright)
   +O_\lambda(Tx).
\end{align}
Next we integrate the error term in~\eqref{phiPart}, which is larger than the error terms of \autoref{thm:lbdSum},
against~\(h_x\mleft(y\mright)\) to obtain
\(
\int_1^\infty O_\lambda(y)h_x(y)\df y = O_{\lambda}(x).\)

Altogether, after dividing \eqref{eqn:lhshxy} and \eqref{eqn:intrhshxy} by \(x\), we have
\begin{align*}
  \sum_{-T < \Im \rho, \, \Im \rho' < T} x^{\rho + \rho' - 1}
  w\mleft( \rho + \rho' - 1 \mright) 
  = \frac{2T\log x}{\pi} \mleft(
  1+\frac{4\cos(\lambda\log x)}{(4+\lambda^2)}
  \mright)
  +O_{f,\lambda,\delta}(T).
\end{align*}
Finally, letting \( x = T^{4 \alpha} \) and dividing both sides
by~\( N_{L_\lambda}\mleft(T\mright)\sim\frac{4}{\pi} T \log T\) yields \autoref{thm:lbdPair}.

\subsection*{Acknowledgements}
The authors thank Machiel van Frankenhuijsen and Jesse Thorner for their
careful reading of a draft of this manuscript and for their insightful comments. The authors also thank an anonymous referee for their diligent reading of the manuscript, whose corrections and suggestions greatly improved the paper.

\printbibliography

\end{document}